\documentclass[12pt]{article}
\usepackage{booktabs}
\usepackage{caption}
\usepackage{mathrsfs}
\usepackage{amsmath}
\usepackage{amsfonts,amsthm,amssymb,mathrsfs,bbding}
\usepackage{graphics,multicol}
\usepackage{graphicx}
\usepackage{color}
\usepackage{enumerate}
\usepackage{caption}
\usepackage{rotating}
\usepackage{lscape}
\usepackage{longtable}
\allowdisplaybreaks[4]
\usepackage{tabularx}
\usepackage[
colorlinks,
anchorcolor=blue,
linkcolor=blue,
citecolor=blue]{hyperref}
\usepackage{extarrows}
\usepackage{cite}
\usepackage{latexsym,bm}
\usepackage{mathtools}
\evensidemargin=\oddsidemargin\topmargin=-1.5cm

\newtheorem{thm}{Theorem}[section]
\newtheorem{prob}{Problem}
\newtheorem{lem}[thm]{Lemma}

\newtheorem{remark}{Remark}

\newtheorem{conj}{Conjecture}

\theoremstyle{definition}

\addtocounter{section}{0}
\begin{document}

\title{\bf Signless Laplacian eigenvalue problems of Nordhaus-Gaddum type}
\author{{Xueyi Huang$^a$ and Huiqiu Lin$^b$\footnote{Corresponding author.}\setcounter{footnote}{-1}\footnote{\emph{E-mail address:} huangxy@zzu.edu.cn (X. Huang), huiqiulin@126.com (H. Lin).}}\\[2mm]
\small $^a$School of Mathematics and Statistics, Zhengzhou University, \\
\small Zhengzhou, Henan 450001, P.R. China\\
\small  $^b$Department of Mathematics, East China University of Science and Technology, \\
\small  Shanghai 200237, P.R. China}

\date{}
\maketitle
{\flushleft\large\bf Abstract} Let $G$ be a graph of order $n$, and let $q_1(G)\geq q_2(G)\geq\cdots\geq q_n(G)$ denote the signless Laplacian eigenvalues of $G$. 
Ashraf and Tayfeh-Rezaie [Electron. J. Combin. 21 (3) (2014)  \#P3.6] showed that $q_1(G)+q_1(\overline{G})\leq 3n-4$, with equality holding if and only if $G$ or $\overline{G}$ is the star $K_{1,n-1}$. In this paper, we discuss the following problem: for $n\geq6$, does $q_2(G)+q_2(\overline{G})\leq 2n-5$ always hold?
We provide positive answers to this problem for the graphs with disconnected complements and the bipartite graphs, and determine the graphs attaining the  bound.  Moreover, we show that $q_2(G)+q_2(\overline{G})\geq n-2$, and the extremal graphs are also characterized.
\begin{flushleft}
\end{flushleft}

\section{Introduction}\label{s-1}

Let $G$ be a simple graph with vertex set $V(G)=\{v_1,\ldots,v_n\}$. We denote the complement  of $G$
by $\overline{G}$, the adjacency matrix of $G$ by $A(G)$, and the  degree (resp. neighborhood) of a vertex $v\in V(G)$  by   $d(v)$ (resp. $N_G(v)$). The \emph{Laplacian matrix} and the \emph{signless Laplacian matrix} of $G$ are the matrices $L(G)=D(G)-A(G)$ and $Q(G)=D(G)+A(G)$,
respectively, where $D(G)=\mathrm{diag}(d(v_1),\ldots,d(v_n))$ is the diagonal  matrix of vertex degrees.
The eigenvalues of $A(G)$, $L(G)$ and $Q(G)$ are called the  \emph{eigenvalues}, \emph{Laplacian eigenvalues}
and \emph{signless Laplacian eigenvalues} ($Q$-\emph{eigenvalues} for short) of $G$,  and  denoted by $\lambda_1(G)\geq\lambda_2(G)\geq\cdots\geq\lambda_n(G)$,
$\mu_1(G)\geq \mu_2(G)\geq\cdots\geq \mu_n(G)$ and $q_1(G)\geq q_2(G)\geq\cdots\geq q_n(G)$, respectively.
Clearly, both $L(G)$ and $Q(G)$ are positive  semidefinite matrices.

Throughout this paper, we denote  the disjoint union of  graphs $G$ and $H$  by $G\cup H$, the disjoint union of $k$'s copies of $G$ by $kG$,   the  join of $G$ and $H$ by $G\nabla H$ which is  obtained from $G\cup H$ by connecting all edges between $G$ and $H$, and the Cartesion product of $G$ and $H$ by $G\square H$. Also, we denote the complete bipartite graph with two parts of sizes $s,t$ by $K_{s,t}$, and the path,   cycle 
and  complete graph on $n$ vertices by $P_n$, $C_n$ and $K_n$, respectively.

In \cite{NG}, Nordhaus and Gaddum considered the lower and upper bounds on the sum and on the product of chromatic number of a graph and its complement. Since then, any bound on the sum or the product of an invariant in a graph $G$ and the same invariant in its complement $\overline{G}$ is called a \emph{Nordhaus-Gaddum type inequality}.
In 2007, Nikiforov \cite{NI} proposed the study of the Nordhaus-Gaddum type
inequalities for all eigenvalues of a graph defining a function given by
$$\max\left\{|\lambda_k(G)|+|\lambda_k(\overline{G})|:|V(G)|=n\right\}, ~\mbox{for}~k=1,\ldots,n.$$
For $k=1$, Nosal \cite{NO} showed that for every graph $G$ of order $n$,
$$\lambda_1(G)+\lambda_1(\overline{G})<\sqrt{2}(n-1).$$
Nikiforov \cite{NI} presented a Nordhaus-Gaddum type result  for the spectral radius of a graph:
$$\frac{4}{3}n-2\le\lambda_1(G)+\lambda_1(\overline{G})<(\sqrt{2}-c)n,$$
where  $c$ is some constant not less than $10^{-7}$. After that, Csikv\'{a}ri \cite{CS}  proved that
$$\lambda_1({G})+\lambda_1(\overline{G})\le\frac{1+\sqrt{3}}{2}n\le1.3661n,$$
which improved the upper bound of Nikiforov. Moreover, Terpai \cite{TE} showed  that
$$\lambda_1(G)+\lambda_1(\overline{G})<\frac{4}{3}n-1.$$
For $k=2$, Nikiforov and Yuan \cite{NY} obtained that $$\lambda_2(G)+\lambda_2(\overline{G})\leq -1+\frac{n}{\sqrt{2}}.$$
Later, Brondani, de Lima and Oliveira \cite{ABLO} posed the following conjecture which improved the bound of Nikiforov and Yuan slightly.
\begin{conj}\label{conj-1}
Let $G$ be a graph on $n$ vertices. Then
$$\lambda_2(G)+\lambda_2(\overline{G})\leq -1+\sqrt{\frac{n^2}{2}-n+1}.$$
\end{conj}

Furthermore, they confirmed  Conjecture \ref{conj-1} for some classes of graphs such as trees, $k$-cyclic graphs, regular bipartite graphs, complete multipartite graphs, generalized line graphs and exceptional graphs.

For the Laplacian eigenvalues, Zhai, Shu and Hong  \cite{ZSH}  (see also You and Liu \cite{YL}) posed the following conjecture on the  Laplacian spread of graphs: 
\begin{conj}\label{conj-2}
Let $G$ be a graph of order $n \geq 2$. Then
$$\mu_1(G)-\mu_{n-1}(G)\leq n-1,$$
with equality holding if and only if $G$ or $\overline{G}$ is isomorphic to the join of an isolated vertex and a disconnected graph of order  $n-1$.
\end{conj}
Notice that the inequality in Conjecture \ref{conj-2} is equivalent to  $\mu_1(G)+\mu_1(\overline{G})\leq 2n-1$ or $\mu_{n-1}(G)+\mu_{n-1}(\overline{G})\geq 1$. Ashraf and Tayfeh-Rezaie \cite{AT} confirmed Conjecture \ref{conj-2} for bipartite graphs, and  Chen and Das \cite{CD} confirmed Conjecture \ref{conj-2} for graphs with $d_1(G)-d_n(G)\leq\sqrt{n-3+2/n}$. Very recently,   Einollahzadeh and Karkhaneei \cite{EK} completely confirmed Conjecture \ref{conj-2}.

With regard to the signless Laplacian eigenvalues,  Ashraf and Tayfeh-Rezaie \cite{AT} showed that $q_1(G)+q_1(\overline{G})\leq 3n-4$, which confirmed a conjecture posed  by Aouchiche and Hansen \cite{AH}. In this paper, we study the similar problem on $q_2(G)+q_2(\overline{G})$ as follows.
\begin{prob}\label{prob1}
Let $G$ be a  graph of order $n$.
What are the upper bound and lower bound of $q_2(G)+q_2(\overline{G})?$
\end{prob}

Note that $Q(G)+Q(\overline{G})=Q(K_n)$. Then by Weyl's inequality (see Lemma \ref{hermitian} below), we have $q_2(G)+q_2(\overline{G})\geq q_3(K_n)=n-2.$
It seems interesting to characterize the extremal graphs attaining this bound. So we first give the following theorem.
\begin{thm}\label{thm-1}
Let $G$ be a graph of order $n\geq 4$. Then  
$$
q_2(G)+q_2(\overline{G})\geq n-2,
$$
with equality holding if and only if $G=K_n$, $nK_1$, $K_{1,n-1}$, $K_{n-1}\cup K_1$, $(2K_1)\nabla K_{n-2}$ or $K_2\cup (n-2)K_1$.
\end{thm}

It is known that $q_2(G)\leq q_2(K_n)=n-2$. Then the first upper bound of $q_2(G)+q_2(\overline{G})$ is as follows.
\begin{thm}\label{thm-2}
Let $G$ be a connected graph of order $n\geq 2$. Then 
$$
q_2(G)+q_2(\overline{G})\leq 2n-4,
$$
with equality holding if and only if $G\in \{K_2,  P_4, C_4\}$.
\end{thm}

\begin{figure}[t]
\begin{center}
\unitlength 2mm 
\linethickness{0.4pt}
\ifx\plotpoint\undefined\newsavebox{\plotpoint}\fi 
\begin{picture}(71,22)(0,2)
\multiput(1,21)(.03370787,-.03370787){89}{\line(0,-1){.03370787}}
\multiput(4,18)(-.03370787,-.03370787){89}{\line(0,-1){.03370787}}
\put(4,18){\line(1,0){3}}
\put(7,18){\line(1,0){5}}
\put(1,21){\circle*{1}}
\put(1,15){\circle*{1}}
\put(4,18){\circle*{1}}
\put(8,18){\circle*{1}}
\put(12,18){\circle*{1}}
\multiput(18,18)(.03370787,.03370787){89}{\line(0,1){.03370787}}
\multiput(21,21)(.03370787,-.03370787){89}{\line(0,-1){.03370787}}
\multiput(24,18)(-.03370787,-.03370787){89}{\line(0,-1){.03370787}}
\multiput(21,15)(-.03370787,.03370787){89}{\line(0,1){.03370787}}
\put(24,18){\line(1,0){4}}
\put(21,21){\circle*{1}}
\put(18,18){\circle*{1}}
\put(21,15){\circle*{1}}
\put(24,18){\circle*{1}}
\put(28,18){\circle*{1}}
\put(34,18){\line(1,0){12}}
\multiput(38,18)(.03370787,.03370787){89}{\line(0,1){.03370787}}
\multiput(41,21)(.03370787,-.03370787){89}{\line(0,-1){.03370787}}
\put(46,18){\line(1,0){2}}
\put(41,21){\circle*{1}}
\put(34,18){\circle*{1}}
\put(38,18){\circle*{1}}
\put(44,18){\circle*{1}}
\put(48,18){\circle*{1}}
\put(54,18){\line(1,0){16}}
\put(54,18){\circle*{1}}
\put(58,18){\circle*{1}}
\put(62,18){\circle*{1}}
\put(66,18){\circle*{1}}
\put(70,18){\circle*{1}}
\multiput(1,9)(.03370787,-.03370787){89}{\line(0,-1){.03370787}}
\multiput(4,6)(-.03370787,-.03370787){89}{\line(0,-1){.03370787}}
\put(1,3){\line(0,1){6}}
\put(4,6){\line(1,0){8}}
\put(1,9){\circle*{1}}
\put(1,3){\circle*{1}}
\put(4,6){\circle*{1}}
\put(8,6){\circle*{1}}
\put(12,6){\circle*{1}}
\put(18,6){\circle*{1}}
\put(21,9){\circle*{1}}
\put(28,9){\circle*{1}}
\put(28,3){\circle*{1}}
\put(21,3){\circle*{1}}
\multiput(18,6)(.03370787,.03370787){89}{\line(0,1){.03370787}}
\put(21,9){\line(1,0){7}}
\put(28,9){\line(0,-1){6}}
\put(28,3){\line(-1,0){7}}
\multiput(21,3)(-.03370787,.03370787){89}{\line(0,1){.03370787}}
\multiput(34,6)(.03370787,.03370787){89}{\line(0,1){.03370787}}
\put(37,9){\line(1,0){7}}
\put(44,9){\line(0,-1){6}}
\put(44,3){\line(0,1){0}}
\put(44,3){\line(-1,0){7}}
\multiput(37,3)(-.03370787,.03370787){89}{\line(0,1){.03370787}}
\put(37,9){\line(0,-1){6}}
\put(34,6){\circle*{1}}
\put(37,9){\circle*{1}}
\put(37,3){\circle*{1}}
\put(44,9){\circle*{1}}
\put(44,3){\circle*{1}}
\multiput(50,6)(.03370787,.03370787){89}{\line(0,1){.03370787}}
\multiput(53,9)(.03370787,-.03370787){89}{\line(0,-1){.03370787}}
\multiput(56,6)(-.03370787,-.03370787){89}{\line(0,-1){.03370787}}
\multiput(53,3)(-.03370787,.03370787){89}{\line(0,1){.03370787}}
\put(56,6){\line(1,0){4}}
\put(53,9){\circle*{1}}
\put(50,6){\circle*{1}}
\put(53,3){\circle*{1}}
\put(56,6){\circle*{1}}
\put(60,6){\circle*{1}}
\put(53,9){\line(0,-1){6}}
\end{picture}
\caption{Connected graphs  of order $5$ with $q_2(G)+q_2(\overline{G})\in (5,6)$.}
\label{fig-0}
\end{center}
\end{figure}
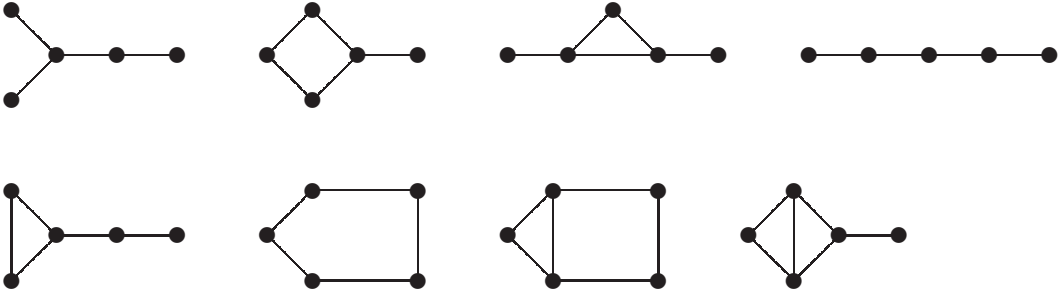

Interesting to us, by using the computer software SageMath v8.7 \cite{ST},  we find that for $n=5$ there are only eight connected graphs (see Figure \ref{fig-0}) with $q_2(G)+q_2(\overline{G})\in (2n-5,2n-4)$, but for $6\leq n\leq 8$ there are no such graphs. Thus we pose the following problem.
\begin{prob}\label{prob2}
Let $G$ be a connected graph of order $n\geq 6$.
Does $q_2(G)+q_2(\overline{G})\leq 2n-5$ always hold?
\end{prob}

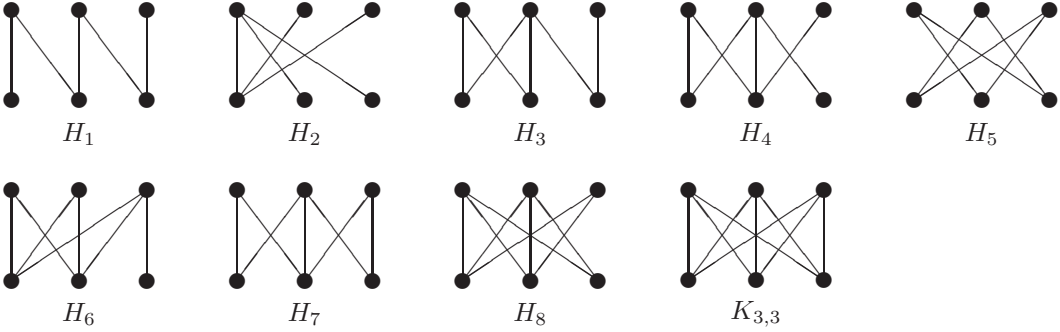
\begin{figure}[t]
\begin{center}
\unitlength 1.5mm 
\linethickness{0.4pt}
\ifx\plotpoint\undefined\newsavebox{\plotpoint}\fi 
\begin{picture}(94,28)(0,-1)
\put(1,26){\circle*{1.5}}
\put(7,26){\circle*{1.5}}
\put(13,26){\circle*{1.5}}
\put(1,18){\circle*{1.5}}
\put(7,18){\circle*{1.5}}
\put(13,18){\circle*{1.5}}
\put(21,26){\circle*{1.5}}
\put(27,26){\circle*{1.5}}
\put(33,26){\circle*{1.5}}
\put(21,18){\circle*{1.5}}
\put(27,18){\circle*{1.5}}
\put(33,18){\circle*{1.5}}
\put(41,26){\circle*{1.5}}
\put(47,26){\circle*{1.5}}
\put(53,26){\circle*{1.5}}
\put(41,18){\circle*{1.5}}
\put(47,18){\circle*{1.5}}
\put(53,18){\circle*{1.5}}
\put(61,26){\circle*{1.5}}
\put(67,26){\circle*{1.5}}
\put(73,26){\circle*{1.5}}
\put(61,18){\circle*{1.5}}
\put(67,18){\circle*{1.5}}
\put(73,18){\circle*{1.5}}
\put(1,18){\line(0,1){8}}
\put(1,26){\line(3,-4){6}}
\put(7,18){\line(0,1){8}}
\put(7,26){\line(3,-4){6}}
\put(13,18){\line(0,1){8}}
\put(21,26){\line(0,-1){8}}
\put(21,18){\line(3,4){6}}
\put(21,18){\line(3,2){12}}
\put(21,26){\line(3,-4){6}}
\put(21,26){\line(3,-2){12}}
\put(41,26){\line(0,-1){8}}
\put(41,18){\line(3,4){6}}
\put(47,26){\line(0,-1){8}}
\put(47,18){\line(-3,4){6}}
\put(47,26){\line(3,-4){6}}
\put(53,18){\line(0,1){8}}
\put(61,26){\line(0,-1){8}}
\put(61,18){\line(3,4){6}}
\put(67,26){\line(0,-1){8}}
\put(67,18){\line(-3,4){6}}
\put(67,26){\line(3,-4){6}}
\put(67,18){\line(3,4){6}}
\put(81,26){\circle*{1.5}}
\put(87,26){\circle*{1.5}}
\put(93,26){\circle*{1.5}}
\put(81,18){\circle*{1.5}}
\put(87,18){\circle*{1.5}}
\put(93,18){\circle*{1.5}}
\put(81,26){\line(3,-4){6}}
\put(87,18){\line(3,4){6}}
\put(93,26){\line(-3,-2){12}}
\put(81,18){\line(3,4){6}}
\put(87,26){\line(3,-4){6}}
\put(93,18){\line(-3,2){12}}
\put(1,10){\circle*{1.5}}
\put(7,10){\circle*{1.5}}
\put(13,10){\circle*{1.5}}
\put(1,2){\circle*{1.5}}
\put(7,2){\circle*{1.5}}
\put(13,2){\circle*{1.5}}
\put(21,10){\circle*{1.5}}
\put(27,10){\circle*{1.5}}
\put(33,10){\circle*{1.5}}
\put(21,2){\circle*{1.5}}
\put(27,2){\circle*{1.5}}
\put(33,2){\circle*{1.5}}
\put(41,10){\circle*{1.5}}
\put(47,10){\circle*{1.5}}
\put(53,10){\circle*{1.5}}
\put(41,2){\circle*{1.5}}
\put(47,2){\circle*{1.5}}
\put(53,2){\circle*{1.5}}
\put(1,10){\line(0,-1){8}}
\put(1,2){\line(3,4){6}}
\put(7,10){\line(0,-1){8}}
\put(7,2){\line(-3,4){6}}
\put(13,10){\line(-3,-2){12}}
\put(13,10){\line(-3,-4){6}}
\put(13,10){\line(0,-1){8}}
\put(27,10){\line(0,-1){8}}
\put(27,2){\line(3,4){6}}
\put(33,10){\line(0,-1){8}}
\put(33,2){\line(-3,4){6}}
\put(27,10){\line(-3,-4){6}}
\put(21,2){\line(0,1){8}}
\put(21,10){\line(3,-4){6}}
\put(41,10){\line(0,-1){8}}
\put(41,2){\line(3,4){6}}
\put(47,10){\line(0,-1){8}}
\put(47,2){\line(-3,4){6}}
\put(41,10){\line(3,-2){12}}
\put(53,2){\line(-3,4){6}}
\put(53,10){\line(-3,-2){12}}
\put(53,10){\line(-3,-4){6}}
\put(61,10){\circle*{1.5}}
\put(67,10){\circle*{1.5}}
\put(73,10){\circle*{1.5}}
\put(61,2){\circle*{1.5}}
\put(67,2){\circle*{1.5}}
\put(73,2){\circle*{1.5}}
\put(61,10){\line(0,-1){8}}
\put(61,2){\line(3,4){6}}
\put(67,10){\line(0,-1){8}}
\put(67,2){\line(-3,4){6}}
\put(61,10){\line(3,-2){12}}
\put(73,2){\line(-3,4){6}}
\put(73,10){\line(-3,-2){12}}
\put(73,10){\line(-3,-4){6}}
\put(73,10){\line(0,-1){8}}
\put(7,15){\makebox(0,0)[cc]{\footnotesize$H_1$}}
\put(27,15){\makebox(0,0)[cc]{\footnotesize$H_2$}}
\put(47,15){\makebox(0,0)[cc]{\footnotesize$H_3$}}
\put(67,15){\makebox(0,0)[cc]{\footnotesize$H_4$}}
\put(87,15){\makebox(0,0)[cc]{\footnotesize$H_5$}}
\put(7,-1){\makebox(0,0)[cc]{\footnotesize$H_6$}}
\put(27,-1){\makebox(0,0)[cc]{\footnotesize$H_7$}}
\put(47,-1){\makebox(0,0)[cc]{\footnotesize$H_8$}}
\put(67,-1){\makebox(0,0)[cc]{\footnotesize$K_{3,3}$}}
\end{picture}
\caption{Connected bipartite graphs satisfying $q_2(G)+q_2(\overline{G})=2n-5$.}
\label{fig-1}
\end{center}
\end{figure}

It is worth mentioning that  for a connected regular graph  $G$ of order $n\geq 6$ the inequality  $q_2(G)+q_2(\overline{G})\leq 2n-5$ always hold. In fact, if $G=K_n$, then we have $q_2(G)+q_2(\overline{G})=n-2<2n-5$. If $G$ is  non-complete  and $k$-regular, Chen and Das \cite{CD} proved that 
$$
\mu_1(G)-\mu_{n-1}(G)< \sqrt{\frac{2nk(n-k-1)}{n-1}},
$$
or equivalently, 
$$
q_2(G)+q_2(\overline{G})< n-2+\sqrt{\frac{2nk(n-k-1)}{n-1}}.
$$
From this inequality, we can deduce that $q_2(G)+q_2(\overline{G})<2n-5$ for $n\geq 9$. For  $6\leq n\leq 8$, by using SageMath v8.7, we find that $q_2(G)+q_2(\overline{G})\leq 2n-5$, and $C_6$, $K_{3,3}$, $K_3\square K_2$ and  $(2K_2)\nabla(3K_1)$ are the only connected regular graphs attaining this bound.

In order to give other positive answers to Problem \ref{prob2}, in this paper, we show the following three results.
\begin{thm}\label{thm-3}
Let $G$ be a connected graph of order $n\geq 6$. If $\overline{G}$ is disconnected, then  
$$
q_2(G)+q_2(\overline{G})\leq 2n-5,
$$
with equality holding if and only if  $G$ is $(K_2\cup K_{n-3})\nabla K_1$, $(2K_2)\nabla(3K_1)$, $K_{3,3}$ or $(K_1\cup K_2)\nabla(K_1\cup K_2)$.
\end{thm}

\begin{thm}\label{thm-4}
Let $G$ be a connected bipartite graph of order $n\geq 6$. Then
$$
q_2(G)+q_2(\overline{G})\leq 2n-5,
$$
with equality holding if and only if $G$ is one of the graphs shown in Figure \ref{fig-1}.
\end{thm}

\begin{thm}\label{thm-5}
Let $G$ be a connected graph of order $n\geq 6$. If $q_2(G)\leq n-3$, then   
$$
q_2(G)+q_2(\overline{G})\leq 2n-5,
$$
with equality holding if and only if $G$ is one of the graphs shown in Figure \ref{fig-1}.
\end{thm}

\begin{remark}
\emph{
According to Theorems \ref{thm-3} and \ref{thm-5}, to resolve Problem \ref{prob1}, it remains to consider the case that $G$ and $\overline{G}$ are connected graphs with $\min\{q_2(G),q_2(\overline{G})\}>n-3$.
}
\end{remark}

\section{Preliminary lemmas}\label{s-2}

Let $A$ be a  Hermitian matrix of order $n$, and let  $\lambda_1(A)\geq \lambda_2(A)\geq \cdots \geq \lambda_n(A)$ denote its   eigenvalues. The following result is well known.

\begin{lem}(Weyl's inequality, \cite{SO}.)\label{hermitian}
Let $A$ and $B$ be Hermitian matrices of order $n$, and let $1 \leq i \leq  n$ and $1  \leq j  \leq n$. Then
$$
\begin{aligned}
&\lambda_i(A) + \lambda_j(B)\leq \lambda_{i+j-n}(A + B), \text{ if } i + j \geq n + 1;\\
&\lambda_i(A) + \lambda_j(B)\geq \lambda_{i+j-1}(A + B), \text{ if } i + j \leq n + 1.
\end{aligned}
$$ 
In either of these inequalities equality holds if and only if there exists a nonzero $n$-vector that is an eigenvector
to each of the three involved eigenvalues.
\end{lem}

Let $\alpha_1\geq \alpha_2\geq \cdots\geq \alpha_n$ and $\beta_1\geq \beta_2\geq \cdots\geq \beta_m$ be two sequences of real numbers  with $m<n$. The second sequence is said to be \textit{interlace} the first one if  
$$\alpha_i\geq \beta_i\geq \alpha_{n-m+i} \text{ for } i=1,2,\ldots,m.$$ 


\begin{lem}(See \cite{HA}.)\label{interlacing-1}
If $B$ is a principal submatrix of a real symmetric matrix $A$, then the eigenvalues of $B$ interlace those of $A$.
 \end{lem}

Let $A$ be a real symmetric matrix of order $n$, and let $X=\{1,2,\ldots,n\}$. For any partition $\Pi:X=X_1\cup \cdots \cup X_m$, the  matrix $A$  can be correspondingly partitioned as
$$
A=\left[
\begin{matrix}
A_{1,1} & \cdots & A_{1,m}\\
\vdots & \ddots & \vdots\\
A_{m,1} & \cdots & A_{m,m}\\
\end{matrix}
\right].
$$
The \emph{characteristic matrix} of  $\Pi$ is the $n\times m$ matrix $\chi_\Pi$  whose columns are the characteristic  vectors of $X_1,\ldots,X_m$, and the \textit{quotient matrix}  of $A$ with respect to $\Pi$ is the matrix $B_\Pi=(b_{ij})_{m\times m}$ with  $b_{ij}=\frac{1}{|X_i|}\mathbf{e}_{|X_i|}^TA_{i,j}\mathbf{e}_{|X_j|}$, where $\mathbf{e}_{|X_i|}$ and $\mathbf{e}_{|X_j|}$ are the all ones $|X_i|$- and $|X_j|$-vectors, respectively.  In particular, the partition $\Pi$ is called \textit{equitable} if each block $A_{i,j}$ has constant row sum.

\begin{lem}(See \cite{BH,GR,HA}.)\label{interlacing-2}
Let  $B_\Pi$ be a quotient matrix of $A$ with respect to some partition $\Pi$. Then  the eigenvalues of $B_\Pi$ interlace those  of $A$. Furthermore, if the partition $\Pi$ is equitable,  then all eigenvalues of $B_\Pi$  are also eigenvalues of $A$, and $A$ has the following two kinds of eigenvectors:
\begin{enumerate}[(i)]
\item the eigenvectors in the column space of $\chi_\Pi$, and the corresponding eigenvalues coincide
with the eigenvalues of $B_\Pi$;

\item the eigenvectors orthogonal to the columns of $\chi_\Pi$, i.e., those eigenvectors that sum to zero on each block $X_i$ for $1\leq i\leq m$.
\end{enumerate}
 \end{lem}
 
\begin{lem}(See \cite{HE}.)\label{interlacing}
Let $G$ be a graph of order $n$, and $H$ a graph obtained from $G$ by deleting some edge. Then
$$q_1(G)\geq q_1(H)\geq q_2(G)\geq q_2(H)\geq\cdots\geq q_{n-1}(G)\geq q_{n-1}(H)\geq q_n(G) \geq q_n(H)\geq 0.$$
\end{lem}

Let $G$ be a graph of order $n$, and let $S=\{v_1,\ldots,v_s\}\subseteq V(G)$ $(s\ge2)$ be a clique (resp. independent set)  such that $N_G(v_i)\setminus S=N_G(v_j)\setminus S$ for all $i,j\in\{1,2,\ldots,s\}$.  Take  $\mathbf{x}_l\in\mathbb{R}^n$ ($2\leq l\leq s$) as the vector defined on $V(G)$ with $\mathbf{x}_l(v_1)=1$, $\mathbf{x}_l(v_l)=-1$ and $\mathbf{x}_l(v)=0$ for $v\not\in\{v_1,v_l\}$, then one can easily verify that $Q(G)\mathbf{x}_l=(d-1)\mathbf{x}_l$ (resp. $Q(G)\mathbf{x}_l=d\mathbf{x}_l$), where $d$ is the common degree of the vertices in $S$. Thus $d-1$ (resp. $d$) is a $Q$-eigenvalue of $G$ with multiplicity at least $s-1$. Then we have

\begin{lem}\label{multiplicity}
If $S$ ($|S|\geq 2$) is a clique (resp. independent set)  of $G$ such that  $N_G(u) \setminus S = N_G(v)\setminus S$ for any $u,v\in S$,  then $G$ has $d-1$ (resp. $d$) as its  $Q$-eigenvalue of  multiplicity at least $|S|-1$, where $d$ is the common degree of the vertices in $S$.
\end{lem}

\begin{lem}(See \cite{FY}.)\label{degree}
Let  $G$ be a connected graph. Then
$$
q_1(G)\leq \max\Bigg\{d(u)+\frac{1}{d(u)}\sum_{uv\in E(G)}d(v):u\in V(G)\Bigg\},
$$
with equality holding if and only if $G$ is either semi-regular bipartite or regular.
\end{lem}
The following lemma can be easily deduced from the Rayleigh's Principle and the Perron-Frobenius Theorem. 
\begin{lem}\label{radius}
Let $G$ be a connected graph of order $n$, and let $u,v$ be two vertices of $G$ which are not adjacent.  Then
$$
q_1(G+uv)>q_1(G).
$$
\end{lem}

A connected bipartite graph is called \textit{balanced} if its vertex classes have the same size,
and \textit{unbalanced} otherwise. Here an isolated vertex can be viewed as an unbalanced bipartite graph with
an empty vertex class.

\begin{lem}(See \cite{LN}.)\label{second}
 If $G$ is a graph of order $n \geq 2$, then 
 $$q_2 (G) \leq  n - 2,$$ 
with equality holding if and only if $\overline{G}$ has a balanced bipartite
component or at least two bipartite components.
\end{lem}

\begin{lem}(See \cite{DAS}.)\label{second-degree}
Let $G$ be a graph with maximum degree $d_1$ and second maximum degree $d_2$. Then
$$q_2(G)\geq d_2-1.$$
If $q_2(G) = d_2-1$, then the maximum and the second maximum degree vertices are adjacent and $d_1=d_2$.
\end{lem}

\begin{lem}(See \cite{GCY}.)\label{least}
If $G$ is a graph on $n$ ($n\geq 6$) vertices and $m$ edges, then
$$q_n(G) \geq \frac{2m}{n-2}-n+1.$$
\end{lem}

\section{Proof of Theorem \ref{thm-1}.}

{\flushleft \bf Proof of Theorem \ref{thm-1}.} 
As $Q(G)+Q(\overline{G})=Q(K_n)$, by Lemma \ref{hermitian}, we have
$$
q_2(G)+q_2(\overline{G})\geq q_3(K_n)=n-2.
$$
It remains to determine the graphs attaining the lower bound. Let $d_1\geq d_2\geq \cdots \geq d_n$ and  $\bar{d}_1\geq \bar{d}_2\geq \cdots \geq \bar{d}_n$ denote the degrees of $G$ and $\overline{G}$, respectively. If $d_2=0$ or $\bar{d}_2=0$, then $G$  is the empty graph $nK_1$ or the complete graph $K_n$, which obviously satisfy $q_2(G)+q_2(\overline{G})=n-2$. Now we always assume that   $d_2,\bar{d}_2\geq 1$ and    $q_2(G)+q_2(\overline{G})=n-2$.  By Lemma \ref{second-degree}, we know that  $q_2(G)\geq d_2-1$ and $q_2(\overline{G})\geq \bar{d}_2-1$.  We consider the  following two cases.

{\flushleft \bf Case 1.} $q_2(G)=d_2-1$ and $q_2(\overline{G})=\bar{d}_2-1$;

In this situation, we have
 $$
 \begin{aligned}
 n-2&=q_2(G)+q_2(\overline{G})\\
 &=d_2-1+\bar{d}_2-1\\
 &=d_2-d_{n-1}+n-3,
 \end{aligned}
 $$
which gives that $d_2=d_{n-1}+1$. Since  $q_2(G)=d_2-1$,  from Lemma \ref{second-degree} we see that $d_1=d_2$, and  all vertices in $G$ of degree $d_2$ must be adjacent to each other. Similarly, by considering the complement  $\overline{G}$, we conclude that $d_{n-1}=d_n$, and the vertices in $G$ of degree $d_{n-1}$ cannot be adjacent to each other. 
Thus $G$ has only two kinds of degrees, i.e., $d_{n-1}+1$ and $d_{n-1}$, and the vertices of degree $d_{n-1}+1$  (resp. $d_{n-1}$) form a clique (resp. independent set). Now partition the vertex set of $G$ as $V(G)=V_1\cup V_2$, where $V_1=\{v\in V(G):d(v)=d_{n-1}+1\}$, $V_2=\{v\in V(G):d(v)=d_{n-1}\}$, and $|V_i|=n_i$ for $i=1,2$.  Since $V_1$ is a clique,  all vertices of $V_1$ have the same number, say $s$, of neighbors in $V_2$. Then we obtain $d_{n-1}+1=n_1-1+s$. Moreover, we claim that $d_{n-1}\leq n_1$ because $V_2$ is an independent set. Thus we get $s\leq 2$. If $s=0$,  there are no edges between $V_1$ and $V_2$,  so $d_{n-1}=0$ and $n_1=2$. This implies that $G$ is just the graph $K_2\cup (n-1)K_1$, which obviously satisfy $q_2(G)=0=d_2-1$, $q_2(\overline{G})=n-2=\bar{d}_2-1$ and $q_2(G)+q_2(\overline{G})=n-2$ by Lemma \ref{second}. If $s=1$, by counting the edges between $V_1$ and $V_2$, we have  $n_1=d_{n-1}n_2$. Combining this with $d_{n-1}+1=n_1-1+s=n_1$, we obtain $d_{n-1}=1$, $n_1=n_2=2$, and so $G=P_4$. However, by simple computation, we have $q_2(P_4)+q_2(\overline{P_4})=2q_2(P_4)=4>2$, a contradiction. If $s=2$, then from $2n_1=d_{n-1}n_2$ and $d_{n-1}+1=n_1-1+s=n_1+1$, we obtain $n_2=2$, $d_{n-1}=n_1=n-2$, and  $G=(2K_1)\nabla K_{n-2}$, which also satisfies the requirement.

{\flushleft \bf Case 2.} $q_2(G)>d_2-1$ or $q_2(\overline{G})>\bar{d}_2-1$.

Without loss of generality, we may assume that $q_2(G)>d_2-1$. We claim that $\bar{d}_2-1\leq q_2(\overline{G})<\bar{d}_2$. In fact, if $q_2(\overline{G})\geq \bar{d}_2$, then  $d_2-d_{n-1}+n-2=d_2-1+\bar{d}_2<q_2(G)+q_2(\overline{G})=n-2$, which gives that $d_2<d_{n-1}$, a contradiction. Also, from $n-2=q_2(G)+q_2(\overline{G})>d_2-1+\bar{d}_2-1=d_2-d_{n-1}+n-3$ we can deduce that $d_2=d_{n-1}$, implying that $G$ has at most three kinds of degrees. Let $V(G)=\{v_1,v_2,\ldots,v_n\}$ with $d(v_1)=d_1$, $d(v_2)=d(v_3)=\cdots=d(v_{n-1})=d_2=d_{n-1}$ and $d(v_n)=d_n$. Then we can partition the vertex set as $V(G)=V_1\cup V_2\cup V_3$, where $V_1=\{v_{1}\}$, $V_2=\{v_2,v_3,\ldots,v_{n-1}\}$ and $V_3=\{v_n\}$.

{\flushleft \bf Subcase 2.1.} $d_2-1<q_2(G)<d_2$;

We claim that $v_1$ is adjacent to all vertices of $V_2$, since otherwise 
$
\left[
\begin{matrix}
d_1& 0\\
0& d_2
\end{matrix}
\right]
$
will be  the principle submatrix of $Q(G)$, which implies that $q_2(G)\geq d_2$ by Lemma \ref{interlacing-1}, a contradiction. Similarly, as $q_2(\overline{G})<\bar{d}_2=n-1-d_{n-1}$, we see that $v_n$ is not adjacent to any vertex of $V_2$.  Observe that  $G[V_1\cup V_2]\cup K_1=(K_1\nabla G[V_2])\cup K_1$ is a spanning subgraph of $G$, and $Q(K_1\nabla G[V_2])$ has the quotient matrix 
$$
B_1=\left[
\begin{matrix}
n-2 & n-2\\
1 & 2d_2-1
\end{matrix}
\right].
$$
Then, by Lemmas \ref{interlacing-2} and  \ref{interlacing}, we have
\begin{equation}\label{equ-4}
q_2(G)\geq \lambda_2(B_1)= \frac{1}{2}\left[n+2d_2-3- \sqrt{n^2 - (4d_2-2)n + 4d_2^2 + 4d_2 - 7}\right].
\end{equation}
Similarly, $\overline{G[V_2\cup V_3]}\cup K_1=(\overline{G[V_2]} \nabla K_1)\cup K_1$ is a spanning subgraph of $\overline{G}$, and $Q(\overline{G[V_2]} \nabla K_1)$ has the quotient matrix
$$
B_2=\left[
\begin{matrix}
n-2 & n-2\\
1 & 2n-2d_2-3
\end{matrix}
\right],
$$
which implies that 
\begin{equation}\label{equ-5}
q_2(\overline{G})\geq \lambda_2(B_2)= \frac{1}{2}\left[3n-2d_2-5- \sqrt{n^2 - (4d_2-2)n + 4d_2^2 + 4d_2 - 7}\right].
\end{equation}
Combining (\ref{equ-4}), (\ref{equ-5})  and $q_2(G)+q_2(\overline{G})=n-2$, we obtain
$$
n-2\geq 2n-4- \sqrt{n^2 - (4d_2-2)n + 4d_2^2 + 4d_2 - 7},
$$
or equivalently,  
\begin{equation}\label{equ-6}
4d_2^2-(4n-4)d_2+6n-11\geq 0.
\end{equation}
Let $f(x)=4x^2-(4n-4)x+6n-11$. For $n\geq 7$, we have $f(2)=f(n-3)=13-2n<0$. Thus  $d_2<2$ or $d_2>n-3$ by (\ref{equ-6}). Since $1\leq d_2\leq n-2$, we have $d_2=1$ or $d_2=n-2$. If $d_2=1$, then $G[V_2]$ is the empty graph,  so $G=K_{1,n-1}$ or $K_{1,n-2}\cup K_1$. In both case, we have $q_2(G)=1$, contrary to $q_2(G)<d_2=1$. If $d_2=n-2$, then $G=K_1\nabla(K_{n-2}\cup K_1)$ or $K_{n-1}\cup K_1$, which are impossible because  $q_2(\overline{K_1\nabla(K_{n-2}\cup K_1)})=q_2(K_{1,n-2}\cup K_1)=1=\bar{d}_2$ and $q_2(K_{n-1}\cup K_1)=n-3=d_2-1$. For $4\leq n\leq 6$, by using  SageMath v8.7, we find that there are no graphs satisfying  $d_2-1<q_2(G)<d_2$, $\bar{d}_2-1\leq q_2(\overline{G})<\bar{d}_2$ and $q_2(G)+q_2(\overline{G})=n-2$.

{\flushleft \bf Subcase 2.2.} $q_2(G)\geq d_2$.

In this situation, we claim that $q_2(G)=d_2$ and $q_2(\overline{G})=\bar{d}_2-1$, since otherwise we have $n-2=q_2(G)+q_2(\overline{G})>d_2+\overline d_2-1=d_2-d_{n-1}+n-2$, a contradiction. As in Subcase 2.1, we see that $v_n$ is not adjacent to any vertex of $V_2$ in $G$.  Note that $\overline{G[V_2\cup V_3]}\cup K_1=(\overline{G[V_2]} \nabla K_1)\cup K_1$ is a spanning subgraph of $\overline{G}$, and $Q(\overline{G[V_2]} \nabla K_1)$ has the quotient matrix
$$
B_3=\left[
\begin{matrix}
n-2 & n-2\\
1 & 2n-2d_2-5+\frac{2s}{n-2}
\end{matrix}
\right],
$$
where $s=d_1$ or $d_1-1$ is the number of edges between $V_1$ and $V_2$ in $G$. By  Lemmas \ref{interlacing-2} and \ref{interlacing}, we have 
\begin{equation}\label{equ-8}
q_2(\overline{G})\geq \lambda_2(B_3)= \frac{(3n-7)(n-2)-(2n-4)d_2+2s-\sqrt{\Delta}}{2n-4},
\end{equation}
where $\Delta=n^4 - ( 4d_2 + 6)n^3 + (4d_2^2 + 28d_2 + 4s + 13)n^2 - (64d_2 + 16d_2^2+(8s+64)d_2+20s+12)n +  16d_2^2+(16s+48)d_2+4s^2+24s+4$. Combining  (\ref{equ-8})  with the fact that $q_2(\overline{G})=\bar{d}_2-1=n-2-d_{n-1}=n-2-d_2$, we obtain
$$
\sqrt{\Delta}\geq n^2-5n+2s+6,
$$
or equivalently, 
$$
(n - 2)d_2^2 - ( n^2 - 5n + 2s + 6)d_2 + n^2 - 4n + 4\geq 0.
$$
Since $s\geq d_1-1\geq d_2-1$, we have
$$
(n - 2)d_2^2 - ( n^2 - 5n + 2(d_2-1) + 6)d_2 + n^2 - 4n + 4\geq 0,
$$
that is,
\begin{equation}\label{equ-9}
(n - 4)d_2^2 -  (n^2 - 5n + 4)d_2 + n^2 - 4n + 4\geq 0.
\end{equation}
Let $g(x)=(n - 4)x^2 -  (n^2 - 5n + 4)x + n^2 - 4n + 4$. For $n\geq 8$, we have  $g(2)=g(n-3)=-(n-5)^2+5<0$, which implies that $d_2<2$ or $d_2>n-3$. Thus $d_2=1$ or $n-2$ due to $1\leq d_2\leq n-2$. If $d_2=1$, then $G=K_{1,s+1}\cup (\frac{n-s-2}{2}K_2)$ with $0\leq s\leq n-2$ or $G=K_{1,s}\cup (\frac{n-s-2}{2}K_2)\cup K_1$ with $1\leq s\leq n-2$. If $s<n-2$, then $G$ has at least one balanced bipartite component $K_2$, implying that $q_2(\overline{G})=n-2$ by Lemma \ref{second},  contrary to the assumption that $q_2(\overline{G})=\bar{d}_2-1=n-3$. Then $s=n-2$, and so $G=K_{1,n-1}$ or $K_{1,n-2}\cup K_1$. We claim that the later case cannot occurs, since otherwise $G$ has two bipartite components, which gives that $q_2(\overline{G})=n-3$ again by Lemma \ref{second}, a contradiction. Thus we have $G=K_{1,n-1}$, which clearly satisfy  the conditions $q_2(G)=d_2=1$ and $q_2(\overline{G})=\bar{d}_2-1=n-3$. If $d_2=n-2$,  then $G[V_2]$ is a clique and $v_1$ is adjacent to all vertices of $V_2$, implying that $G=K_{1}\nabla(K_{n-2}\cup K_1)$ or $K_{n-1}\cup K_1$. In both cases, we can deduce a contradiction because $q_2(\overline{K_{1}\nabla(K_{n-2}\cup K_1)})=q_2(K_1 \cup K_{1,n-2})=1>0=\bar{d}_2-1$ and $q_2(K_{n-1}\cup K_1)=n-3<n-2=d_2$. For $4\leq n\leq 8$, by using SageMath v8.7, we find that $K_{1,n-1}$ is the only graph satisfying $q_2(G)=d_2$, $q_2(\overline{G})=\bar{d}_2-1$  and $q_2(G)+q_2(\overline{G})=n-2$.

Therefore, we conclude that $G=K_{1,n-1}$ under the assumption that $q_2(G)>d_2-1$. Also, by considering  $\overline{G}$ with $q_2(\overline{G})>\bar{d}_2-1$,  we obtain $G=K_{n-1}\cup K_1$. 

Concluding the above results, we have $q_2(G)+q_2(\overline{G})=n-2$ if and only if $G$ is $K_n$, $nK_1$, $K_{1,n-1}$, $K_{n-1}\cup K_1$, $(2K_1)\nabla K_{n-2}$ or  $K_2\cup (n-2)K_1$.

We complete the proof.\qed

\section{Proof of Theorems \ref{thm-2}--\ref{thm-5}.}


{\flushleft \bf Proof of Theorem \ref{thm-2}.} 
Notice that  $Q(G)+Q(\overline{G})=Q(K_n)$. By Lemma \ref{hermitian}, we have
\begin{equation*}
\left\{
\begin{aligned}
&q_2(G)+q_n(\overline{G})\leq q_2(K_n)\\
&q_2(\overline{G})+q_n(G)\leq q_2(K_n)
\end{aligned}
\right.,
\end{equation*}
which implies that
\begin{equation}\label{equ-1}
\left\{
\begin{aligned}
&q_2(G)\leq n-2-q_n(\overline{G})\leq n-2\\
&q_2(\overline{G})\leq n-2-q_n(G)\leq n-2
\end{aligned}
\right.
\end{equation}
because $q_2(K_n)=n-2$ and $q_n(G), q_n(\overline{G})\geq 0$. Summing up the two inequalities  in (\ref{equ-1}) side by side, we obtain 
\begin{equation*}
q_2(G)+q_2(\overline{G})\leq 2n-4,
\end{equation*}
where the  equality holds if and only if $q_2(G)=q_2(\overline{G})=n-2$ (and thus $q_n(G)=q_2(\overline{G})=0$). By Lemma \ref{second}, this is the case that both  $G$ and $\overline{G}$ contain  a balanced bipartite component or at least two bipartite components. Then $G$ must be a connected balanced bipartite graph of order $n=2s$ for some positive integer $s$. We claim that $s\leq 2$, since otherwise $\overline{G}$  contains $2K_s$ as its spanning subgraph, implying that $q_n(\overline{G})\geq q_n(2K_s)=s-2\geq 1$, contrary to $q_n(\overline{G})=0$. Thus $G=K_2$, $P_4$ or $C_4$. Conversely, one can check that $K_2$, $P_4$ and $C_4$ satisfy the relation $q_2(G)+q_2(\overline{G})= 2n-4$. 

We compete the proof. \qed


{\flushleft \bf Proof of Theorem \ref{thm-3}.} 
By the assumption, we  suppose that  $\overline{G}=H_1\cup H_2\cup \cdots \cup H_k$, where  $H_i$ is  connected and $n_i=|V(H_i)|$ for $1\leq i\leq k$. Assume that  $n_1\geq n_2\geq \cdots\geq n_k$. We consider the following two situations.

{\flushleft\bf Case 1.} $k\geq 3$; 

If $n_1=1$, then $\overline{G}$ is the empty graph,  so $G=K_n$ and $q_2(G)+q_2(\overline{G})=n-2<2n-5$. Now suppose $n_1\geq 2$. Since $n_1\leq n-(k-1)\leq n-2$ and $n_2\leq \frac{n-1}{2}$,  we have
\begin{equation}\label{equ-2}
\begin{aligned}
q_2(\overline{G})&\leq \max\{q_2(H_1),q_1(H_2),q_1(H_3),\ldots,q_1(H_k)\}\\
&\leq \max\{n_1-2,q_1(K_{n_2})=2n_2-2\}\\
&\leq \max\{n-4,n-3\}\\
&=n-3,
\end{aligned}
\end{equation}
where the second inequality follows from Lemma \ref{interlacing}. Moreover, we have $q_2(G) \leq n-2$  by (\ref{equ-1}), and thus
$$
q_2(G)+q_2(\overline{G})\leq 2n-5.
$$
Here the equality holds if and only if $q_2(\overline{G})=n-3$ and $q_2(G)=n-2$. If $q_2(\overline{G})=n-3$, from (\ref{equ-2}) we know that $n_2=\frac{n-1}{2}$  ($n$ is odd). Since $n_1\geq n_2$ and $n_3\geq 1$, we have $k=3$,  $n_1=n_2=\frac{n-1}{2}$ and $n_3=1$. Also, we see that $q_2(\overline{G})=q_1(H_2)=q_1(K_{n_2})$, and so $H_2=K_{n_2}=K_{\frac{n-1}{2}}$ by Lemma \ref{radius}. Moreover, we have $n-3=q_2(\overline{G})\leq q_1(\overline{G})=q_1(H_1)\leq q_1(K_{\frac{n-1}{2}})=n-3$, which implies that $H_1=K_{\frac{n-1}{2}}$ again by Lemma \ref{radius}. Therefore, we have $\overline{G}=(2K_{\frac{n-1}{2}}) \cup K_1$. By Lemma  \ref{second}, this gives that $q_2(G)<n-2$ because $n\geq 6$ is odd. Hence, we conclude that $q_2(G)+q_2(\overline{G})<2n-5$  in this situation.

{\flushleft\bf Case 2.} $k=2$. 

In this situation, we have $3\leq n_1\leq n-1$ and  $n_2\leq \frac{n}{2}$.  We distinguish  the following two subcases to discuss.

{\flushleft\bf Subcase 2.1.} $n_2\leq \frac{n-1}{2}$; 

We have 
\begin{equation}\label{equ-3}
\begin{aligned}
q_2(\overline{G})&\leq \max\{q_2(H_1),q_1(H_2)\}\\
&\leq \max\{n_1-2,q_1(K_{n_2})=2n_2-2\}\\
&\leq \max\{n-3,n-3\}\\
&=n-3.
\end{aligned}
\end{equation}
As above, we have $ q_2(G)+q_2(\overline{G})\leq 2n-5$, where the equality holds if and only if $q_2(\overline{G})=n-3$ and $q_2(G)=n-2$. 

Now suppose that $q_2(\overline{G})=n-3$ and $q_2(G)=n-2$. First assume that  $q_2(H_1)\geq q_1(H_2)$.  According to (\ref{equ-3}), we have $n-3=q_2(\overline{G})=q_2(H_1)\leq n_1-2\leq n-3$, which implies that $n_1=n-1$, $q_2(H_1)=n-3$ and $H_2=K_1$. As $q_2(G)=n-2$,  by Lemma \ref{second}, the graph $H_1$ must be bipartite. Also, since $H_1$ has  $n-1$ vertices and $q_2(H_1)=n-3$, we see that $\overline{H_1}$ has a  balanced bipartite component or at least two bipartite components. If $H_1$ is a complete bipartite graph, say $H_1=K_{s,n-1-s}$ ($1\leq s\leq n-2$ and $n\geq 6$), then $\overline{H_1}=K_s\cup K_{n-1-s}$ cannot have two bipartite components, and so must have a balanced bipartite component, that is, $\overline{H_1}=K_2\cup K_{n-3}$.  Thus we have  $\overline{G}=K_{2,n-3}\cup K_1$, i.e., $G=(K_2\cup K_{n-3})\nabla K_1$, which obviously satisfy the conditions $q_2(\overline{G})=n-3$ and $q_2(G)=n-2$. If $H_1$ is not complete bipartite, then $\overline{H_1}$ is connected, and so must be a  balanced bipartite graph by the above arguments, which is impossible due to $n-1\geq 5$.  Now assume that  $q_2(H_1)<q_1(H_2)$. Then we have $n-3=q_2(\overline{G})=q_1(H_2)\leq q_1(K_{n_2})=2n_2-2\leq n-3$, implying that $n_2=\frac{n-1}{2}$ ($n$ must be odd), $n_1=\frac{n+1}{2}$, and $H_2=K_{\frac{n-1}{2}}$ by Lemma \ref{radius}. Also, since $q_2(G)=n-2$, we claim that $H_1$ is a balanced bipartite graph by Lemma \ref{second}. Furthermore, we have $q_1(\overline{G})=q_1(H_1)$ because   $q_2(\overline{G})=q_1(H_2)$. Then $n-3=q_2(\overline{G})\leq q_1(\overline{G})=q_1(H_1)\leq q_1(K_{\frac{n_1}{2},\frac{n_1}{2}})=q_1(K_{\frac{n+1}{4},\frac{n+1}{4}})=\frac{n+1}{2}$, implying that  $n=7$ because $n\geq 6$ is odd. For $n=7$, since both sides of the above inequality equal to $4$, we must have   $H_1=K_{2,2}=C_4$ by  Lemma \ref{radius}. Thus $\overline{G}=K_{2,2}\cup K_3$, i.e., $G=(2K_2)\nabla (3K_1)$, which obviously satisfy the relation $q_2(\overline{G})=4$ and $q_2(G)=5$. 

Therefore, in this subcase, we conclude that  $q_2(G)+q_2(\overline{G})\leq 2n-5$, where the equality holds if and only if $G=(K_2\cup K_{n-3})\nabla K_1$ or $(2K_2)\nabla (3K_1)$. 

{\flushleft\bf Subcase 2.2.} $\frac{n-1}{2}< n_2\leq \frac{n}{2}$.

We see that $n$ must be even and $n_2=\frac{n}{2}=n_1$.  If $q_2(\overline{G})=q_2(H_1)$ or $q_2(H_2)$, then $q_2(\overline{G})\leq \frac{n}{2}-2<n-3$, and so $q_2(G)+q_2(\overline{G})<2n-5$. Thus we may assume that $q_2(\overline{G})=\min\{q_1(H_1),q_1(H_2)\}$. If $q_2(\overline{G})<n-3$, there is nothing to prove. We only need to consider the situation that $q_2(\overline{G})\geq n-3$, from which we obtain  $q_1(H_i)\geq n-3$ for $i=1,2$. Let $d_1$ denote the maximum degree of $H_i$. By Lemma \ref{degree}, we have $q_1(H_i)\leq 2d_1$, which gives that $d_1\geq \frac{n-3}{2}$. Thus $d_1=\frac{n}{2}-1$ due to $d_1\leq \frac{n}{2}-1$ and $n$ is even. For $i=1,2$, let $w_i\in V(H_i)$ be such that
$$
d(w_i) +\frac{1}{d(w_i)}\sum_{w_iv\in E(H_i)}d(v)=\max\Bigg\{d(u)+\frac{1}{d(u)}\sum_{uv\in E(H_i)}d(v):u\in V(H_i)\Bigg\}.
$$
We claim that $d(w_i)=d_1$ or $d_1-1$, since otherwise we get $n-3\leq q_1(H_i)\leq d(w_i) +\frac{1}{d(w_i)}\sum_{w_iv\in E(H_i)}d(v)\leq d_1-2+d_1=n-4$ by Lemma \ref{degree}, a contradiction.  

If $d(w_i)=d_1=\frac{n}{2}-1$,  we obtain $\sum_{w_iv\in E(H_i)}d(v)\geq (\frac{n}{2}-2)(\frac{n}{2}-1)$ by using the fact that $n-3\leq q_1(H_i)\leq d(w_i) +\frac{1}{d(w_i)}\sum_{w_iv\in E(H_i)}d(v)$. Thus $2m(H_i)\geq d(w_i) +\sum_{w_iv\in E(H_i)}d(v)\geq(\frac{n}{2}-1)^2$. If $d(w_i)=d_1=\frac{n}{2}-2$, as above, we have $\sum_{w_iv\in E(H_i)}d(v)\geq (\frac{n}{2}-2)(\frac{n}{2}-1)$. Then, for each $v\in N_{H_i}(w_i)$, we have $d(v)=\frac{n}{2}-1$. Let $w_i'$ be the unique vertex that is not adjacent to $w_i$ in $H_i$. We see that $d(w_i')=\frac{n}{2}-2$ because  all vertices of $N_{H_i}(w_i)$ are adjacent to $w_i'$. Thus we have $H_i=K_{\frac{n}{2}-2}\nabla(2K_1)$, and  $2m(H_i)=d(w_i')+d(w_i)+\sum_{w_iv\in E(H_i)}d(v)=(\frac{n}{2}-2)(\frac{n}{2}+1)$. For $n\geq 12$, in both cases we can deduce that 
$$
q_{\frac{n}{2}}(H_i)\geq  \frac{2m(H_i)}{\frac{n}{2}-2}-\frac{n}{2}+1>1
$$
by  Lemma \ref{least}. Thus  $q_n(\overline{G})=\min\{q_{\frac{n}{2}}(H_1),q_{\frac{n}{2}}(H_2)\}>1$, and we have $q_2(G)<n-3$ according to (\ref{equ-1}). Hence, we conclude that $q_2(G)+q_2(\overline{G})<2n-5$ because $q_2(\overline{G})\leq n-2$. For each even $n$ with $6\leq n\leq 11$,  by  considering the condition $q_1(H_i)\geq n-3$ and using   SageMath v8.7,  we obtain $H_i=K_{\frac{n}{2}}$ or  $K_{\frac{n}{2}-2}\nabla(2K_1)$  for $i=1,2$. Recall that $\overline{G}=H_1\cup H_2$. By simple computation, we find that $q_2(G)+q_2(\overline{G})\leq 2n-5$, and the equality holds if and only if $\overline{G}=2K_3$ or $2(K_1\nabla (2K_1))=2K_{1,2}$, i.e., $G=K_{3,3}$ or $(K_2\cup K_1)\nabla (K_2\cup K_1)$.

Therefore, in this subcase, we conclude that  $q_2(G)+q_2(\overline{G})\leq 2n-5$, with equality holding if and only if $G=K_{3,3}$ or $(K_2\cup K_1)\nabla (K_2\cup K_1)$.

We complete the proof.
\qed

\begin{figure}[t]
\begin{center}
\unitlength 1.5mm 
\linethickness{0.4pt}
\ifx\plotpoint\undefined\newsavebox{\plotpoint}\fi 
\begin{picture}(50,18)(0,-1)
\put(2,4){\circle*{1.5}}
\put(12,4){\circle*{1.5}}
\put(20,4){\circle*{1.5}}
\put(30,4){\circle*{1.5}}
\put(16,17){\circle*{1.5}}
\put(38,4){\circle*{1.5}}
\put(48,4){\circle*{1.5}}
\multiput(16,17)(-.0362694301,-.0336787565){386}{\line(-1,0){.0362694301}}
\multiput(16,17)(-.033613445,-.109243697){119}{\line(0,-1){.109243697}}
\multiput(16,17)(.033613445,-.109243697){119}{\line(0,-1){.109243697}}
\multiput(16,17)(.0362694301,-.0336787565){386}{\line(1,0){.0362694301}}
\put(34,17){\circle*{1.5}}
\multiput(34,17)(-.033613445,-.109243697){119}{\line(0,-1){.109243697}}
\multiput(34,17)(-.0362694301,-.0336787565){386}{\line(-1,0){.0362694301}}
\multiput(34,17)(.033613445,-.109243697){119}{\line(0,-1){.109243697}}
\multiput(34,17)(.0362694301,-.0336787565){386}{\line(1,0){.0362694301}}
\put(25,4){\oval(16,6)[]}
\put(43,4){\oval(16,6)[]}
\put(7,4){\oval(16,6)[]}
\put(6,4){\circle*{1.5}}
\put(24,4){\circle*{1.5}}
\put(42,4){\circle*{1.5}}
\multiput(6,4)(.0336700337,.0437710438){297}{\line(0,1){.0437710438}}
\multiput(16,17)(.033613445,-.054621849){238}{\line(0,-1){.054621849}}
\multiput(24,4)(.0336700337,.0437710438){297}{\line(0,1){.0437710438}}
\multiput(34,17)(.033613445,-.054621849){238}{\line(0,-1){.054621849}}
\put(9,4){\makebox(0,0)[cc]{$\cdots$}}
\put(27,4){\makebox(0,0)[cc]{$\cdots$}}
\put(45,4){\makebox(0,0)[cc]{$\cdots$}}
\put(16,19){\makebox(0,0)[cc]{\footnotesize$u$}}
\put(34,19){\makebox(0,0)[cc]{\footnotesize$v$}}
\put(2,2){\makebox(0,0)[cc]{\footnotesize$u_1$}}
\put(6,2){\makebox(0,0)[cc]{\footnotesize$u_2$}}
\put(12,2){\makebox(0,0)[cc]{\footnotesize$u_{s_1}$}}
\put(20,2){\makebox(0,0)[cc]{\footnotesize$w_1$}}
\put(24,2){\makebox(0,0)[cc]{\footnotesize$w_2$}}
\put(30,2){\makebox(0,0)[cc]{\footnotesize$w_{s_0}$}}
\put(38,2){\makebox(0,0)[cc]{\footnotesize$v_1$}}
\put(42,2){\makebox(0,0)[cc]{\footnotesize$v_2$}}
\put(48,2){\makebox(0,0)[cc]{\footnotesize$v_{s_2}$}}
\put(25,-1){\makebox(0,0)[cc]{\footnotesize$S_{0}$}}
\put(7,-1){\makebox(0,0)[cc]{\footnotesize$S_{1}$}}
\put(43,-1){\makebox(0,0)[cc]{\footnotesize$S_{2}$}}
\end{picture}
\caption{The graph $H(s_{0},s_{1},s_{2})$.}
\label{fig-2}
\end{center}
\end{figure}
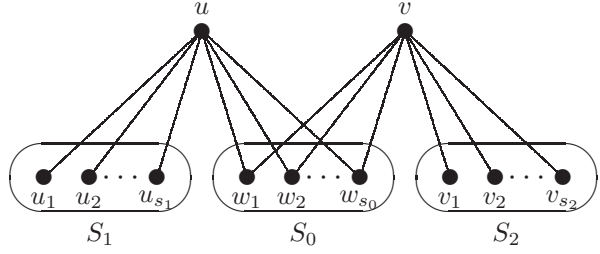


{\flushleft \bf Proof of Theorem \ref{thm-4}.} 
If $\overline{G}$ is disconnected, by Theorem \ref{thm-3}, we have $q_2(G)+q_2(\overline{G})\leq 2n-5$, where the equality holds if and only if $G=K_{3,3}$.  Now suppose that $\overline{G}$ is connected. Let $V(G)=V_1\cup V_2$ be the bipartition of $G$ with $|V_i|=n_i$ for $n=1,2$ and $n_1\geq n_2$. 

If $n_1,n_2\geq 3$, then $q_n(\overline{G})\geq q_n(K_{n_1}\cup K_{n_2})=n_2-2\geq 1$, and so $q_2(G)\leq n-3$ by (\ref{equ-1}). Thus we have $q_2(G)+q_2(\overline{G})\leq 2n-5$, where the the equality holds if and only if $q_2(G)=n-3$ and $q_2(\overline{G})=n-2$.  Now assume that  $q_2(G)=n-3$ and $q_2(\overline{G})=n-2$. By Lemma \ref{second}, the bipartite graph $G$ must be  balanced, i.e., $n_1=n_2=\frac{n}{2}$. Furthermore, since $n-3 = q_2(G)\leq q_2(K_{\frac{n}{2},\frac{n}{2}})=\frac{n}{2}$ and $n\geq 6$, we have $n=6$ and $q_2(G)=3$. By using SageMath v8.7, we find that all connected balanced bipartite graphs of order $6$ with connected complements satisfying $q_2(G)=3$  and $q_2(\overline{G})=4$ are exactly the graphs $H_1$--$H_7$ shown in Figure \ref{fig-1}.  

Next suppose  $n_2\leq 2$. If $n_2=1$, then $G=K_{1,n-1}$, contrary to the connectedness of  $\overline{G}$. Thus $n_2=2$ and $n_1=n-2$. Suppose  $V_2=\{u,v\}$ and $V_1=S_0\cup S_1 \cup S_2$, where $S_{0}=N_G(u)\cap N_G(v)$, $S_{1}=N_G(u)\setminus N_G(v)$, $S_{2}=N_G(v)\setminus N_G(u)$ and $|S_{i}|=s_{i}$ for $i=0,1,2$. Then we see that  $G$ is of  the form $H(s_{0},s_{1},s_{2})$ ($s_{0}+s_{1}+s_{2}=n-2$) shown in Figure \ref{fig-2}. By the connectedness of $G$ and $\overline{G}$, we see that $s_0\geq 1$ and $\max\{s_1,s_2\}\geq 1$. By the symmetry, we only need to consider the following two situation.

{\flushleft\bf Case 1.} $s_1,s_2\geq 1$;

For $6\leq n\leq 8$, by using SageMath v8.7, we find that  $q_2(G)+q_2(\overline{G})<2n-5$ holds for all such $G$'s. Now suppose  $n\geq 9$.   If $s_{1},s_{2}\geq 2$,  then $\overline{G}$ has $2K_3\cup K_{n-6}$ as its spanning subgraph, and so $q_n(\overline{G})\geq q_n(2K_3\cup K_{n-6}) \geq 1$ by Lemma \ref{interlacing}. This implies that $q_2(G)\leq n-3$ according to (\ref{equ-1}). Furthermore, as the connected bipartite graph $G$ is not balanced,  we have $q_2(\overline{G})<n-2$ by Lemma \ref{second}, and therefore,   $q_2(G)+q_2(\overline{G})<2n-5$. It remains to consider the case that $s_{1}=1$ or $s_{2}=1$. Without loss of generality, we may assume that $s_{1}=1$. If $s_{2}\geq 2$, then  $G$ is a spanning subgraph of  $H(n-5,1,2)$. By Lemma \ref{multiplicity}, it is easy to see that  the graph $H(n-5,1,2)$ has  $Q$-eigenvalues  $2$  of multiplicity at least $n-6$ and $1$ of multiplicity at least one. Furthermore, by Lemma \ref{interlacing-2} and the proof of Lemma \ref{multiplicity}, we claim that  the remaining five $Q$-eigenvalues, denoted by $\alpha_1>\alpha_2\geq \alpha_3\geq \alpha_4>\alpha_5=0$,  must be in the   quotient matrix 
$$
B_1=
\left[\begin{matrix}
2 &0& 0& 1& 1\\
 0& 1& 0& 1& 0\\
 0& 0& 1& 0& 1\\
  n-5& 1& 0& n-4& 0\\
 n-5 &0 &2 &0 &n-3
\end{matrix}\right]
\begin{matrix}
S_0\\
S_1\\
S_2\\
u\\
v
\end{matrix}.
$$
By simple computation, the characteristic polynomial of $B_1$ is  $\phi(B_1,x)=xf(x)$, where $f(x)=x^4 - (2n-3)x^3 + (n^2 - n - 4)x^2 - (2n^2-8n +2)x + n^2 - 5n$. Since $f(n-3)=-(n-5)(n-6)<0$, we have $\alpha_2<n-3$ or $\alpha_3>n-3$. We claim that the later case cannot occurs, since otherwise we obtain $2n-3=\mathrm{trace}(B_1)=\alpha_1+\alpha_2+\alpha_3+\alpha_4>3(n-3)$, which is impossible due to $n\geq 9$. Thus $\alpha_2<n-3$, and we have $q_2(G)\leq q_2(H(n-5,1,2))=\max\{\alpha_2,2\}<n-3$. Therefore, we conclude that $q_2(G)+q_2(\overline{G})<2n-5$. If $s_{2}=1$, then $G=H(n-4,1,1)$. As above, we see that $G$ has  $Q$-eigenvalue $2$ of multiplicity at least $n-5$, and the remaining five $Q$-eigenvalues, denoted by $\beta_1>\beta_2\geq \beta_3\geq \beta_4>\beta_5=0$, are contained in the quotient matrix
$$
B_2=
\left[\begin{matrix}
2 &0& 0& 1& 1\\
 0& 1& 0& 1& 0\\
 0& 0& 1& 0& 1\\
  n-4& 1& 0& n-3& 0\\
 n-4 &0 &1 &0 &n-3
\end{matrix}\right]
\begin{matrix}
S_0\\
S_1\\
S_2\\
u\\
v
\end{matrix}.
$$
By simple computation, we get $\beta_2=\frac{n-2+\sqrt{n^2-8n+20}}{2}$, and so $q_2(G)=\frac{n-2+\sqrt{n^2-8n+20}}{2}$. Now consider the complement graph $\overline{G}=\overline{H(n-4,1,1)}$. We see that $\overline{G}$ has $Q$-eigenvalue $n-4$ of multiplicity at least $n-5$, and the remaining  five $Q$-eigenvalues  lie in the quotient matrix
$$
B_3=
\left[\begin{matrix}
2n-8 &1 &1 &0 &0\\
n-4 &n-2 &1 &0 &1\\
 n-4 &1 &n-2 &1 &0\\
  0 &0 &1 &2 &1\\
  0 &1 &0& 1& 2
\end{matrix}\right]
\begin{matrix}
S_0\\
S_1\\
S_2\\
u\\
v
\end{matrix}.
$$
By simple computation, the characteristic polynomial of $B_3$ is  $\phi(B_3,x)=f_1(x)f_2(x)$, where $f_1(x)= x^3 - (3n-6)x^2 + (2n^2 - 3n - 12)x - 6n^2 + 38n - 56$ and $f_2(x)=x^2 -(n - 2)x + n - 4$. Let $\gamma_1\geq \gamma_2\geq \gamma_3$ be the three roots of $f_1(x)$, and $\gamma_1'\geq \gamma_2'$ the two roots of $f_2(x)$. We have $\gamma_1>2n-6$ because  $f_1(2n-6)=-4n+16<0$, and   $\gamma_1'=\frac{n-2+\sqrt{n^2-8n+20}}{2}$, $\gamma_2'=\frac{n-2-\sqrt{n^2-8n+20}}{2}$. Also, since $f_1(\gamma_1')=-2(n - 4)(n - 3-\sqrt{n^2 - 8n + 20})<0$, we obtain $\gamma_1>\gamma_1'>\gamma_2$ or $\gamma_1'<\gamma_3$. We claim that the later case cannot occurs, since otherwise we have $4n-8=\mathrm{trace}(B_3)>\gamma_1+3\gamma_1'+\gamma_2'>4n-10+\sqrt{n^2-8n+20}$, a contradiction. Therefore, we have  $q_2(\overline{G})=\gamma_1'=\frac{n-2+\sqrt{n^2-8n+20}}{2}$, and so $q_2(G)+q_2(\overline{G})=n-2+\sqrt{n^2-8n+20}<2n-5$. 

{\flushleft\bf Case 2.} $s_1=0$ and $s_2\geq 1$.

If $s_2\geq 2$, then $G$ is a spanning subgraph of $H(n-4,0,2)$. We see that $H(n-4,0,2)$ has $Q$-eigenvalues  $2$  of multiplicity at least $n-5$, $1$ of multiplicity at least one, and  the remaining four $Q$-eigenvalues, denoted by $\alpha_1>\alpha_2\geq \alpha_3> \alpha_4=0$,  are in the  quotient matrix 
$$
B_1=
\left[\begin{matrix}
2 & 0&1 &1\\
0& 1&0&1\\
n-4& 0& n-4&0\\
n-4& 2 &0 &n-2
\end{matrix}\right]
\begin{matrix}
S_0\\
S_2\\
u\\
v
\end{matrix}.
$$
By simple computation, the characteristic polynomial of $B_1$ is  $\phi(B_1,x)=xf(x)$, where $f(x)=x^3 - (2n-3)x^2 + (n^2 - 2n - 2)x - n^2 + 4n$. Similarly, as $f(n-3)=-n+6\leq 0$,  we have $\alpha_2\leq n-3$ or $\alpha_3\geq n-3$. We claim that the later case cannot occurs, since otherwise we have $2n-3=\mathrm{trace}(B_1)=\alpha_1+\alpha_2+\alpha_3+\alpha_4>3\alpha_3\geq 3(n-3)$, a contradiction. Thus  $q_2(G)\leq q_2(H(n-4,0,2))=\alpha_2\leq n-3$.  Furthermore, by Lemma \ref{second}, we have  $q_2(\overline{G})<n-2$  because $G$ is connected but not balanced. Thus $q_2(G)+q_2(\overline{G})<2n-5$. If $s_2=1$, then $G=H(n-3,0,1)$. For $n=6,7$, one can easily check that $q_2(G)+q_2(\overline{G})<2n-5$. Now suppose $n\geq 8$. We see that $G$ has $Q$-eigenvalue $2$ of multiplicity at least $n-4$, and the  remaining four $Q$-eigenvalues, denoted by $\beta_1>\beta_2\geq \beta_3>\beta_4=0$, are contained in the quotient matrix 
$$
B_2=
\left[\begin{matrix}
2 & 0&1 &1\\
0& 1&0&1\\
n-3& 0& n-3&0\\
n-3& 1 &0 &n-2
\end{matrix}\right]
\begin{matrix}
S_0\\
S_2\\
u\\
v
\end{matrix}.
$$
The characteristic polynomial of $B_2$ is  $\phi(B_2,x)=xg(x)$, where $g(x)=x^3  -(2n-2)x^2 + (n^2 - n - 2)x - n^2 + 3n$. Since $g(n-\frac{5}{2})=\frac{-2n+15}{8}<0$, as above, we have $q_2(G)=\beta_2<n-\frac{5}{2}$.  Now consider the complement graph $\overline{G}=\overline{H(n-3,0,1)}$. We see that $\overline{G}$ has $Q$-eigenvalue $n-4$ of multiplicity at least $n-4$, and the remaining  four $Q$-eigenvalues, denoted by $\gamma_1>\gamma_2\geq \gamma_3 \geq \gamma_4>0$  lie in the quotient matrix
$$
B_3=
\left[\begin{matrix}
2n-7 & 1&0 &0\\
n-3& n-2&1&0\\
0& 1& 2&1\\
0& 0 &1 &1
\end{matrix}\right]
\begin{matrix}
S_0\\
S_2\\
u\\
v
\end{matrix}.
$$
The characteristic polynomial of $B_3$ is $\phi(B_3,x)=x^4 -(3n-6)x^3 + (2n^2 - 3n - 10)x^2 - (6n^2 - 35n + 48)x + 2n^2 - 14n + 24$. Also, we have $\phi(B_3,2n-6)=-4(n-3)(n-4)<0$, and $\phi(B_3,n-\frac{5}{2})=-\frac{1}{16}((2n - 11)(4n^2 - 24n + 39))<0$, which implies that $\gamma_1>2n-6$ and $\gamma_2<n-\frac{5}{2}$. Thus $q_2(\overline{G})= \max\{n-4,\gamma_2\}<n-\frac{5}{2}$. Hence, we get $q_2(G)+q_2(\overline{G})<2n-5$.

Concluding the above results, we obtain that $q_2(G)+q_2(\overline{G})\leq 2n-5$, where the equality holds if and only if $G$ is one of the graphs shown in Figure \ref{fig-1}.

We complete the proof.\qed


{\flushleft \bf Proof of Theorem \ref{thm-5}.} 
Clearly, we have  $q_2(G)+q_2(\overline{G})\leq 2n-5$, where the equality holds if and only if $q_2(G)=n-3$ and $q_2(\overline{G})=n-2$. If $q_2(\overline{G})=n-2$, by Lemma \ref{second}, the connected graph $G$ must be bipartite and balanced. Thus the results follows from the proof of Theorem \ref{thm-4} immediately.
\qed




{\flushleft \bf Acknowledgements.} The first author is  supported  by the National Natural Science Foundation of China (No. 11671344), the China Postdoctoral Science Foundation and the Postdoctoral Research Sponsorship in Henan Province (No. 1902011). The second author is supported by the National Natural Science Foundation of China (Nos. 11771141) and Fundamental Research Fund for the Central Universities (No. 222201714049).

\end{document}